\documentclass[10pt,psamsfonts]{amsart}
\usepackage{amsmath}
\usepackage{amsthm}
\usepackage{amssymb}
\usepackage{amscd}
\usepackage{amsfonts}
\usepackage{amsbsy}
\usepackage{graphicx}
\usepackage{url}
\usepackage{overpic}
\usepackage{hyperref}
\usepackage{float}
\usepackage{xcolor}
\usepackage[pagewise]{lineno}
\usepackage[foot]{amsaddr}

\newcommand{\R}{\ensuremath{\mathbb{R}}}
\newcommand{\N}{\ensuremath{\mathbb{N}}}

\newcommand{\ov}{\overline}

\newcommand{\T}{\theta}

\newcommand{\bP}{\mathbf{P}}
\newcommand{\x}{\mathbf{x}}

\def\e{\varepsilon}

\newtheorem {theorem} {Theorem}

\newtheorem {corollary}{Corollary}

\newtheorem {example} {Example}
\newtheorem {remark}{Remark}

\newtheorem {mtheorem} {Theorem}

\newtheorem {problem} {Problem}

\newcommand\blfootnote[1]{%
	\begingroup
	\renewcommand\thefootnote{}\footnote{#1}%
	\addtocounter{footnote}{-1}%
	\endgroup
}

\usepackage{mathpazo}

\begin{document}

\title[A version of Hilbert's 16th Problem for 3D polynomial vector fields]
{A version of Hilbert's 16th Problem\\ for 3D polynomial vector fields:\\
counting isolated invariant tori}

\author[D.D. Novaes and P.C.C.R. Pereira]
{Douglas D. Novaes and Pedro C.C.R. Pereira}

\address{Universidade Estadual de Campinas (UNICAMP), Departamento de Matem\'{a}tica, Instituto de Matemática, Estatística e Computação Científica (IMECC) - Rua S\'{e}rgio Buarque de Holanda, 651, Cidade
Universit\'{a}ria Zeferino Vaz, 13083--859, Campinas, SP, Brazil}
\email{ddnovaes@unicamp.br}
\email{pedro.pereira@ime.unicamp.br}



\begin{abstract}
 Hilbert's 16th Problem, about the maximum number of limit cycles of planar polynomial vector fields of a given degree $m$, has been one of the most important driving forces for new developments in the qualitative theory of vector fields. Increasing the dimension, one cannot expect the existence of a finite upper bound for the number of limit cycles of, for instance, $3$D polynomial vector fields of a given degree $m$. Here, as an extension of such a problem in the $3$D space, we investigate the number of isolated invariant tori in $3$D polynomial vector fields. In this context, given a natural number $m$, we denote by $N(m)$ the upper bound for the number of isolated invariant tori of $3$D polynomial vector fields of degree $m$. Based on a recently developed averaging method for detecting invariant tori, our first main result provides a mechanism for constructing $3$D differential vector fields with a  number $H$ of normally hyperbolic invariant tori from a given planar differential vector field with $H$ hyperbolic limit cycles. The strength of our mechanism in studying the number $N(m)$ lies in the fact that the constructed $3$D differential vector field is polynomial provided that the given planar differential vector field is polynomial. Accordingly, our second main result establishes a lower bound for $N(m)$ in terms of lower bounds for the number of hyperbolic limit cycles of planar polynomial vector fields of degree $[m/2]-1$.  Based on this last result, we apply a methodology due to Christopher \& Lloyd to show that $N(m)$ grows as fast as $m^3/128$. Finally, the above-mentioned problem is also formulated for higher dimensional polynomial vector fields.
\end{abstract}
\blfootnote{
	DDN is partially supported by S\~{a}o Paulo Research Foundation (FAPESP) grants 2022/09633-5, 2019/10269-3, and 2018/13481-0, and by Conselho Nacional de Desenvolvimento Cient\'{i}fico e Tecnol\'{o}gico (CNPq) grants 438975/2018-9 and 309110/2021-1. PCCRP is supported by S\~{a}o Paulo Research Foundation (FAPESP) grant 2020/14232-4. \\
	\rule{60pt}{0.4pt}}
\maketitle


\section{Introduction}

Given a natural number $m$, the second part of Hilbert's 16th Problem inquires about the existence of a uniform upper bound $H(m)$ (named {\it Hilbert Number}) for the amount of limit cycles of planar polynomial vector fields of degree $m$. Formally, let $\pi(P,Q)$ denote the number of limit cycles of the polynomial differential system $\dot x=P(x,y)$, $\dot y=Q(x,y)$, then,
\[
H(m):=\sup\{\pi(P,Q):\textrm{degree}(P),\textrm{degree}(Q)\leq m\}.
\]

The only known Hilbert Number is for $m=1$, because planar linear vector fields do not admit limit cycles, that is $H(1)=0$. For $m\geq2$, it is still unknown even if $H(m)$ is finite. Since Hilbert's list publication \cite{Hilbert}, this problem has been one of the most important driving forces for new developments in the qualitative theory of vector fields (see, for instance, \cite{Ilyashenko02,L88}). 

Simpler versions of Hilbert’s 16th Problem have been proposed. One worth mentioning is Smale’s 13th Problem \cite{smale00} (see, also, \cite{LMP77}), which focuses on finding an upper bound for the number of limit cycles in Liénard polynomial differential equations, i.e. $P(x, y) = y - F(x)$ and $Q(x, y) = -x$, where $F(x)$ is a polynomial of degree $m$.

\subsection{Lower bounds for the Hilbert Numbers}
So far, the main results obtained focus exclusively on lower bounds for $H(m)$, even for the simplified version of Hilbert’s 16th Problem proposed by Smale (see, for instance, \cite{MD11,MH15,DPR07}). In \cite{CL95}, Christopher \& Lloyd showed that $H(m)$ grows as fast as $n^2\log n$. This result has been revisited by several authors (see, for instance, \cite{lvarez2020,HL12,LCC02}). In \cite{HL12}, Han \& Li  improved Christopher \& Lloyd's result by showing that $H(m)$ grows as fast as $(m+2)^2\log (m+2)/(2\log 2)$, which means that
\begin{equation}\label{liminf}
\lim_{m\to\infty}\inf\dfrac{H(m)}{(m+2)^2\log(m+2)}\geq\dfrac{1}{2\log 2}
\end{equation}
For $m\in\N$ small, the current best known lower bounds for $H(m)$ were provided by Prohens \& Torregrosa in \cite{PT19}. For the reader's convenience, we collect some of them here: $H(2)\geq 4$, $H(3)\geq13$, $H(4)\geq 28, H(5)\geq37, H(6)\geq53, H(7)\geq74, H(8)\geq96, H(9)\geq120,$ $H(10)\geq142$, $H(13)\geq212, H(17)\geq384, H(21)\geq568, H(31)\geq1184, H(35)\geq1536, H(39)\geq1920,$ and $H(43)\geq2272$.

\begin{remark}\label{rem}
The asymptotic estimate \eqref{liminf} and the lower bounds above can be realized by hyperbolic limit cycles.  Indeed, Han \& Li \cite{HL12}   detected two types of limit cycles, namely, small-amplitude limit cycles and large-amplitude limit cycles. The small-amplitude limit cycles are obtained from degenerate Hopf-bifurcations and, therefore, they are hyperbolic (see \cite[Lemma 2.5]{HL12}). In particular, the asymptotic estimate \eqref{liminf} is obtained only by considering small-amplitude limit cycles (see \cite[Theorem 3.4]{HL12}). The limit cycles obtained by Prohens \& Torregrosa in \cite{PT19} are also generated from degenerate Hopf-bifurcations and, therefore, they are hyperbolic.
\end{remark}

In light of Remark \ref{rem}, let $\pi_h(P,Q)$ denote the number of hyperbolic limit cycles of the polynomial differential system $\dot x=P(x,y)$, $\dot y=Q(x,y)$. Then, for a natural number $m$, define
\[
H_h(m):=\sup\{\pi_h(P,Q):\textrm{degree}(P),\textrm{degree}(Q)\leq m\}.
\]
We also denote by $\underline{H}_h(m)$ a finite lower bound for $H_h(m)$. In the asymptotic estimate \eqref{liminf} and in the lower bounds above, $H$ can be changed by $\underline{H}_h$.

\subsection{A version of Hilbert's 16th Problem in the $3$D space}\label{problem3d}
Increasing the dimension, one cannot expect the existence of a finite upper bound for the number of limit cycles of, for instance, $3$D polynomial vector fields of a given degree $m$. Indeed, for $m=2$, the R\"{o}ssler system \cite{Rossler76},
\begin{equation}\label{eq:rs}
\left\{\begin{aligned}
\dot x&=-y-x,\\
\dot y&=x+a\, y,\\
\dot z&=b+z\,(x-c),
\end{aligned}\right.
\end{equation}
undergoes cascades of period-doubling bifurcations due to the existence of  Shilnikov homoclinic connections for some values of the parameters (see, for instance, \cite{algaba,BBS}). This implies the number of limit cycles becomes unbounded by varying the parameters.

Here, as an extension of Hilbert's 16th Problem in the $3$D space, we suggest to increase also the dimension of the sought invariant object by inquiring about the number of isolated invariant tori in $3$D polynomial vector fields. In this context, we denote by $\tau(P,Q,R)$ the number of isolated invariant tori of the polynomial differential system $\dot x=P(x,y,z)$, $\dot y=Q(x,y,z)$, $\dot z=R(x,y,z)$. Then, for a natural number $m$, we define
\[
N(m):=\sup\{\tau(P,Q,R):\textrm{degree}(P),\textrm{degree}(Q),\textrm{degree}(R)\leq m\}.
\]
We also define
\[
N_h(m):=\sup\{\tau_h(P,Q,R):\textrm{degree}(P),\textrm{degree}(Q),\textrm{degree}(R)\leq m\},
\]
where $\tau_h(P,Q)$ denotes the number of normally hyperbolic invariant tori of the polynomial differential system $\dot x=P(x,y,z)$, $\dot y=Q(x,y,z)$, $\dot z=R(x,y,z)$. 
As with the Hilbert Number, we do not know whether $N(m)$ and $N_h(m)$ are finite or not.

\section{Statements of the main results}

\subsection{$3$D differential systems with normally hyperbolic invariant tori}\label{sec:prel}

Our first main result provides a mechanism for constructing $3$D differential systems with a  number $H$ of normally hyperbolic invariant tori from a given planar differential system with $H$ hyperbolic limit cycles. The strength of our mechanism in studying the number $N(m)$ lies in the fact that the constructed $3$D differential system is polynomial provided that the given planar differential system is polynomial. So far, to the best of our knowledge, there was not a simple way of constructing $3$D polynomial differential systems with a given number of isolated invariant tori. 

\begin{mtheorem}\label{thm:mainA}
Let $P,Q:\R^2\to\R$ be smooth functions. Assume that the planar differential system
\begin{equation}\label{eq:2dpde}
\begin{cases}
\dot x=P(x,y),\\
\dot y=Q(x,y),
\end{cases}
\end{equation}
has at least $H$ hyperbolic limit cycles contained in the region $K=(0,b)\times(\alpha,\beta)$, for some $b>0$ and $\beta>\alpha$. Consider the following one parameter family of $3$D differential systems:
\begin{equation}\label{eq:3dpde}
\begin{cases}
\dot x=-y,\\
\dot y=x+\e\,y\, P(x^2+y^2,z),\\
\dot z=2\,\e\,y^2 Q(x^2+y^2,z).
\end{cases}
\end{equation}
Then, there exists $\bar \e>0$ such that the differential system \eqref{eq:3dpde} has at least $H$ normally hyperbolic invariant tori for every $\e\in(0,\bar\e)$. Moreover, if $\gamma=\{(u(t),v(t)),t\in[0,T)\}\subset K$ is a parametrization of a limit cycle of \eqref{eq:2dpde}, then for each $\e\in(0,\bar\e)$ there exists a normally hyperbolic invariant torus $T_{\e}^{\gamma}$ of  \eqref{eq:3dpde} with the same stability as $\gamma$ and converging, as $\e$ goes to $0$, to the torus $T_0^{\gamma}:=\{(u(t)\cos(\T),u(t)\sin(\T),v(t)):(t,\T)\in[0,T)\times [0,2\pi)\}$ which is invariant for the unperturbed system \eqref{eq:3dpde}$_{\e=0}$.
\end{mtheorem}

The proof of Theorem \ref{thm:mainA} is provided in Subsection \ref{sec:proofA} and it is based on a recently developed  averaging method for detecting normally hyperbolic invariant tori \cite{Pereira2023}  which is discussed in Subsection \ref{sec:at}.

\begin{remark}\label{rem:trans}
The technical assumption that the hyperbolic limit cycles are contained in the region $K$ is not restrictive. Indeed, since we are dealing only with a finite number $H$ of limit cycles, then this condition can always be assumed without loss of generality just by considering a translation of the variable $x$.
\end{remark}

In the sequel, we present an example illustrating the application of Theorem \ref{thm:mainA} in the context of Liénard polynomial systems.
\begin{example}
Consider a planar differential system of Liénard type, given by $\dot{u} = v - F(u)$ and $\dot{v} = -u$. Through the transformation $(x, y) = (a+v, u )$, this Liénard system is equivalent to system \eqref{eq:2dpde} with $P(x, y) = -y$ and $Q(x, y) = x - a - F(y)$. As a result, the corresponding $3$D system \eqref{eq:3dpde} becomes
\begin{equation}\label{eq:3dpdelienard}
\begin{cases}
\dot{x} = -y,\\
\dot{y} = x - \e\,y\,z,\\
\dot{z} = 2\,\e\,y^2 \left(x^2 + y^2 - a - F(z)\right),
\end{cases}
\end{equation}
which has degree equal to $\max\{4, m+2\}$, provided that $F$ is a polynomial of degree $m$.

We claim that, for each $m \geq 6$, there exist a polynomial function $F$ of degree $m$ and parameters $a > 0$ and $\e > 0$ such that the $3$D polynomial differential system \eqref{eq:3dpdelienard} has $m - 2$ normally hyperbolic invariant tori.

Indeed, by \cite[Theorem 1]{MH15}, for each $m \geq 6$, there exists a polynomial $F$ of degree $m$ such that the corresponding Liénard system and, consequently, the planar system $\dot{x} = -y$ and $\dot{y} = x - a - F(y)$ has $m - 2$ hyperbolic limit cycles. Moreover, the parameter $a > 0$ can be chosen so that these limit cycles are contained within the half-plane $x > 0$. Thus, by applying Theorem \ref{thm:mainA}, the claim follows.

With this example, we can further conclude that, for $m \geq 6$, $N(m + 2) \geq m - 2$, which in turn implies $N(m) \geq m - 4$ for $m \geq 8$. In what follows, we will improve this lower bound.
\end{example}

\subsection{Lower bounds for $N(m)$}
Our second main result  provides a lower bound for $N_h(m),$ and consequently also for $N(m)$, in terms of lower bounds for the number of hyperbolic limit cycles of planar polynomial vector fields of degree $[m/2]-1$.

\begin{mtheorem}\label{thm:mainB}
The relationship
\begin{equation}\label{eq:thm}
N_h(m)\geq \underline{H}_h\left(\left[\dfrac{m}{2}\right]-1\right)
\end{equation}
holds for each natural number $m\geq 2$.
\end{mtheorem}

Theorem \ref{thm:mainB} allows us to provide asymptotic estimates for $N_h$ based on the asymptotic estimates for $\underline{H}_h$. Indeed, consider the increasing function $F(m)=(m+2)^2\log (m+2)$. Since $[m/2]\geq (m-1)/2$, Theorem \ref{thm:mainB} implies that
\[
\dfrac{N_h(m)}{F((m-1)/2-1)}\geq \dfrac{\underline{H}_h\left(\left[\dfrac{m}{2}\right]-1\right)}{F([m/2]-1)}.
\]
Thus, taking \eqref{liminf} and Remark \ref{rem} into account, we obtain the following result.
\begin{corollary}\label{cor}
$N_h(m)$ grows as fast as
\[
\dfrac{1}{8\log 2}(m+1)^2\log\left(\dfrac{m+1}{2}\right).
\]
\end{corollary}

In Section \ref{sec:fc}, we improve Corollary \ref{cor} by showing that $N_h(m)$ actually grows at least as fast as $m^3/128$.

For $m\in\N$ small, better lower bounds for $N_h(m)$ can be given based on better lower bounds for $\underline{H}_h$. We collect some of them in the next result.

\begin{corollary} \label{cor2}
$N_h(7)\geq N_h(6)\geq \underline{H}_h(2)\geq 4$, 
$N_h(9)\geq N_h(8)\geq \underline{H}_h(3)\geq 13$,
$N_h(11)\geq N_h(10)\geq \underline{H}_h(4)\geq 28$,
$N_h(13)\geq N_h(12)\geq \underline{H}_h(5)\geq 37$,
$N_h(15)\geq N_h(14)\geq \underline{H}_h(6)\geq 53$,
$N_h(17)\geq N_h(16)\geq \underline{H}_h(7)\geq 74$,
$N_h(19)\geq N_h(18)\geq \underline{H}_h(8)\geq 89$,
$N_h(21)\geq N_h(20)\geq \underline{H}_h(9)\geq 120$, 
$N_h(23)\geq N_h(22)\geq \underline{H}_h(10)\geq 142$,
$N_h(29)\geq N_h(28)\geq \underline{H}_h(13)\geq 212$,
$N_h(37)\geq N_h(36)\geq \underline{H}_h(17)\geq 348$,
$N_h(45)\geq N_h(44)\geq \underline{H}_h(21)\geq 568$,
$N_h(65)\geq N_h(64)\geq \underline{H}_h(31)\geq 1184$,
$N_h(73)\geq N_h(72)\geq \underline{H}_h(35)\geq 1536$,
$N_h(81)\geq N_h(80)\geq \underline{H}_h(39)\geq 1920$, and
$N_h(89)\geq N_h(88)\geq \underline{H}_h(43)\geq 2272$.
\end{corollary}

As we shall see in the proof of Theorem \ref{thm:mainB}, the lower bound \eqref{eq:thm} relies on Theorem \ref{thm:mainA}, which will provide a $3$D polynomial vector field of degree $2k+2$ with $\underline{H}_h(k)$ normally hyperbolic invariant tori from a given planar vector field of degree $k$ with $\underline{H}_h(k)$ hyperbolic limit cycles. Thus, it is natural to seek for lower bounds of $N$ in terms of $H$ better than the one provided by Theorem \ref{thm:mainB}. This leads us to formulate the following problem:

\begin{problem}
Given a natural number $m$, does there exist a $3$D polynomial vector field with degree $m$ having more isolated invariant tori than $H\left(\left[\dfrac{m}{2}\right]-1\right)$. If so, can these invariant tori be taken normally hyperbolic?
\end{problem}

\section{Proof of the main results}
The proof of Theorem \ref{thm:mainA} is based on an averaging method for finding normally hyperbolic invariant tori, which was recently developed in \cite{Pereira2023} (see also \cite{CN20,Novaes2023}). We start this section discussing such method. 
Theorems \ref{thm:mainA} and \ref{thm:mainB} are, then, proved in Subsections  \ref{sec:proofA} and \ref{sec:proof}, respectively.

\subsection{Averaging method for finding invariant tori}\label{sec:at}

Consider the following $T$-periodic non-autonomous differential equations given in the following standard form:
\begin{align} \label{eq:e1}
    \dot \x = \varepsilon F(t, \x,\e) \quad (t,\x,\e)\in \R \times D \times [0,\varepsilon_0],
\end{align}
where $D$ is an open bounded subset of $\R^2,$  $\varepsilon_0>0,$ and $F:\R \times D \times [0,\varepsilon_0]\to\R^2$ is a $C^r$ function, with $r\geq 2$, $T$-periodic in the variable $t$. 

The next result provides sufficient conditions for the existence of normally hyperbolic  invariant tori in the extended phase space of the differential equation \eqref{eq:e1}, that is $\mathbb{S}^1\times D$, with $\mathbb{S}^1\equiv\R/(T\mathbb{Z})$. 

\begin{theorem} \label{thm:nhit}
Consider the differential equation \eqref{eq:e1}. Suppose that the guiding system 
\begin{equation}\label{guiding}
\dot \x=\dfrac{1}{T}\int_0^T F_1(t,\x)d\x
\end{equation}
has an attracting hyperbolic limit cycle $\gamma$. Then, there exists $\ov \e>0$ such that, for each $\e\in(0,\ov \e)$, the differential equation \eqref{eq:e1} has a $T$-periodic solution $\gamma_{int}$ and a normally hyperbolic attracting invariant torus in the extended phase space. In addition, the torus surrounds the periodic solution $\gamma_{int}$ and converges to $\mathbb{S}^1\times\gamma$ as $\e$ goes to $0$.
\end{theorem}

A weaker version of Theorem \ref{thm:nhit}, without the conclusion of the normal hyperbolicity of the invariant torus, follows from a result due to Hale (see \cite[Theorem 18-2]{halebook}). 

\subsection{Proof of Theorem \ref{thm:mainA}}\label{sec:proofA}
Applying the cylindrical change of coordinates $(x,y,z)=(r\cos(\T),r\sin(\T),z)$, the differential system \eqref{eq:3dpde} becomes
\begin{equation}\label{eq:3dpde2}
\begin{cases}
\dot \T=1+\e\,\sin(\T)\cos(\T)P(r^2,z),\\
\dot r=\e\, r\sin^2(\T) P(r^2,z),\\
\dot z=2\,\e\, r^2\sin^2(\T) Q(r^2,z).
\end{cases}
\end{equation}
Notice that, for the open bounded set $D=(0,\sqrt{b})\times (\alpha,\beta)$, $b>0$, there exists $\e_0>0$ such that $\dot \T>0$ for every $(\T,r,z)\subset [0,2\pi]\times \ov D $ and $\e\in[0,\e_0]$. Thus, we can take $\T$ as the new independent variable, by applying a time-rescale that divides the differential equations in \eqref{eq:3dpde2} by $\dot\T$. By doing so, one gets the following $2\pi$-periodic non-autonomous differential system
\begin{equation}\label{eq:napde}
\begin{cases}
\dfrac{d r}{d\T}=\dfrac{ \e\,r\sin^2(\T) P(r^2,z)}{1+\e\,\sin(\T)\cos(\T)P(r^2,z)}=\e\, r \sin^2(\T)P(r^2,z)+\e^2R_1(\T,r,z,\e),\vspace{0.1cm}\\
\dfrac{ d z}{d \T}=\dfrac{2\,\e\, r^2\sin^2(\T) Q(r^2,z)}{1+\e\,\sin(\T)\cos(\T)P(r^2,z)}=\e\,2r^2 \sin^2(\T) Q(r^2,z)+\e^2R_2(\T,r,z,\e).
\end{cases}
\end{equation}
Notice that \eqref{eq:napde} is written in the standard form for applying the averaging method, to wit
\[
\dot x=\e F_1(\T,\x)+\e^2 R(\T,\x,\e), \,\, (\T,\x,\e)\in \R \times D \times[0,\e_{0}],
\] 
where 
\[
\begin{aligned}
&T=2\pi,\quad t=\T,\quad \x=(r,z),\\ &F_1(\T,\x)=\Big( r \sin^2(\T)P(r^2,z),2r^2 \sin^2(\T) Q(r^2,z)\Big),\text{ and}\\
 &R(\T,\x,\e)=\Big(R_1(\T,r,z,\e),R_2(\T,r,z,\e)\Big).
\end{aligned}
\] 
Taking \eqref{guiding} into account, the guiding system is given by
\begin{equation}\label{eq:gs}
(r',z')=\dfrac{1}{2\pi}\int_0^{2\pi} F_1(\T,\x)d\T=\left(\dfrac{1}{2}r P(r^2,z),r^2Q(r^2,z)\right),\quad (r,z)\in D.
\end{equation}

In what follows, we show that \eqref{eq:gs} has at least $H$ hyperbolic limit cycles. To do that, consider the transformation $(\rho,z)=(r^2,z),$ which is a diffeomorphism on $\{(r,z):r>0\}$. Applying this transformation to the guiding system \eqref{eq:gs}, we get
\begin{equation}\label{eq:gs2}
(\rho',z')=\left(\rho P(r,z),\rho Q(\rho,z)\right), \quad (\rho,z)\in K=(0,b)\times(\alpha,\beta).
\end{equation}
Finally, proceeding with a time-rescaling, dividing the right-hand side of \eqref{eq:gs2} by $\rho$, we conclude that the differential system \eqref{eq:gs2} and, therefore, the guiding system \eqref{eq:gs} is $C^{\infty}$-equivalent to the planar polynomial differential system \eqref{eq:2dpde}, which has $H$ hyperbolic limit cycles.

The proof concludes by applying Theorem \ref{thm:nhit} and going back through the changes of variables.

\subsection{Proof of Theorem \ref{thm:mainB}}\label{sec:proof}

Given a natural number $k$, by the definition of  $\underline{H}_h(k)$,  there exist polynomials of degree $k$, $P(x,y)$ and $Q(x,y)$, for which the planar polynomial differential system \eqref{eq:2dpde} has at least $\underline{H}_h(k)$ hyperbolic limit cycles. As pointed out in Remark \ref{rem:trans}, by considering a translation in the variable $x$, we can assume, without loss of generality, that all the $\underline{H}_h(k)$ limit cycles of \eqref{eq:2dpde} lie in a region $K=(0,b)\times(\alpha,\beta)$, for some $b>0$ and $\beta>\alpha$.

From Theorem \ref{thm:mainA}, there exists $\ov\e>0$ such that the $3$D differential system \eqref{eq:3dpde} has at least $H$ normally hyperbolic invariant tori for every $\e\in(0,\bar\e)$. Moreover, \eqref{eq:3dpde} is polynomial of degree $2k+2$. Hence, $N(2k+2)\geq \underline{H}_h(k)$. 

Thus, by taking $k=n-1$, one has $N(2n)\geq\underline{H}_h(n-1)$, which implies \eqref{eq:thm} for $m=2n$ even. 

Now, for $m$ odd, let us say $m=2n+1,$ one has 
\[N(m)=N(2n+1)\geq N(2n)\geq\underline{H}_h(n-1).\]
Since, in this case, $[m/2]=n,$  then the relationship \eqref{eq:thm} also holds for $m$ odd.

\section{Further comments}

\subsection{Improving the asymptotic growth of $N(m)$}\label{sec:fc}
The methodology introduced by Christopher \& Lloyd \cite{CL95} can be adapted to improve the asymptotic growth estimation of $N_h(m)$ given by Corollary \ref{cor}. Indeed, from Corollary \ref{cor2}, there exists a $3$D polynomial vector field $X_0=(X_{01},X_{02},X_{03})$ of degree $6$ with $4$ normally hyperbolic invariant tori. By performing a translation, we can assume, without loss of generality, that such tori are located in the set $S=(-1,+\infty)\times(-1,+\infty)\times(-1,+\infty)$. Considering the following transformation 
\[
\phi:(x,y,z)\mapsto (x^2-1,y^2-1,z^2-1),
\]
define 
\[
\begin{aligned}
X_1(x,y,z)=&2\,x\,y\,z\,(D\phi(x,y,z))^{-1} X_0\circ\phi(x,y,z)\\
=&\left(y\,z\,X_{01}\circ\phi(x,y,z)\,,\,x\,z\,X_{02}\circ\phi(x,y,z)\,,\,x\,y\,X_{03}\circ\phi(x,y,z)\right).
\end{aligned}
\]
Notice that $X_1$ is a polynomial vector field of degree $m_1=2\cdot 6+2=14$, which corresponds to the pullback of $X_0$ by $\phi$ multiplied by the factor $2x\,y\,z$, that represents a reparametrization of time in each octant. Of course, $\phi$ is not a diffeomorphism on the whole $\R^2$. Nevertheless, it is indeed a diffeomorphism between each octant of $\R^3$ onto $S$. This means that $X_1$ restricted to each octant is equivalent to $X_0$ restricted to $S$. Therefore, $X_1$ has at least $\tau_1=8\cdot 4=32$ normally hyperbolic invariant tori, $4$ of them in each one of the $8$ octants. This procedure can be continued in order to create a sequence of polynomial vector fields $(X_k)$ of degree $m_k$ and with at least $\tau_k$ normally hyperbolic invariant tori, where $m_k=2 m_{k-1}+2$ and $\tau_k=8 \tau_{k-1}$. Since $m_0=6$ and $\tau_0=4$, it yields
\[
m_k=2(2^{k+2}-1)\quad\text{and}\quad \tau_k=2^{2+3k}=\dfrac{(m_k+2)^3}{128}.
\]
This means that
\[
N_h(m_k)\geq \dfrac{(m_k+2)^3}{128}>\frac{m_k^3}{128},
\]
which implies that 
\begin{equation}\label{liminf2}
\lim_{m\to\infty}\inf\dfrac{N_h(m)}{m^3}\geq\dfrac{1}{128},
\end{equation}
that is, $N_h(m)$ grows as fast as $m^3/128$.

\subsection{Into higher dimensions}
The problem introduced in Section \ref{problem3d} can be extended to higher dimensions. Let $\bP=(P_1,\ldots,P_{d+1}):\R^{d+1}\to\R^{d+1}$, where each $P_i$ is a polynomial in the variable $\x$. Define the degree of $\bP$ as $\textrm{degree}(\bP)=\max_i \textrm{degree}(P_i)$.

We denote by $\tau^d(\bP)$ (resp. $\tau_h^d(\bP)$) the number of isolated (resp. normally hyperbolic) invariant $d$-dimensional tori ($d$-tori) of the polynomial differential system $\dot \x=\bP(\x)$. Then, for a natural number $m$, we define:
\[
N^d(m):=\sup\{\tau^d(\bP):\textrm{degree}(\bP)\leq m\}
\]
and
\[
N_h^d(m):=\sup\{\tau_h^d(\bP):\textrm{degree}(\bP)\leq m\}.
\]
Notice that $N^1(m)=H(m)$, $N^1_h(m)=H_h(m)$, $N^2(m)=N(m)$, and $N^2_h(m)=N_h(m)$.

Based on the previous discussion, we expect that $N^d_h(m)$ grows at least as fast as $L\,m^{d+1}$ for some $L>0$. Indeed, once we have an example of a $(d+1)$-dimensional polynomial vector field $X_0$ of degree $m_0$ with at least $\tau_0$ normally hyperbolic invariant $d$-tori, we can proceed as follows:

First, by performing a translation, we can assume, without loss of generality, that such tori are located in the set $S=(-1,+\infty)^{d+1}$. Consider the following transformation $(x_1,\ldots,x_{d+1})\mapsto (x_1^2-1,\ldots,x_{d+1}^2-1)$ and define 
\[
\begin{aligned}
X_1(x_1,&\ldots,x_{d+1})\\
=&2 x_1\cdots x_{d+1}(D\phi(x_1,\ldots,x_{d+1}))^{-1} X_0\circ\phi(x_1,\ldots,x_{d+1})\\
=&\left(x_2\ldots x_{d+1}X_{01}\circ\phi(x_1,\ldots,x_{d+1}),\,\ldots\, , x_1\ldots x_{d}\,X_{0(d+1)}\circ\phi(x_1,\ldots,x_{d+1})\right).
\end{aligned}
\]
Notice that $X_1$ is a polynomial vector field of degree $m_1=2\cdot m_0+d$, corresponding to the pullback of $X_0$ by $\phi$ multiplied by the factor $2 x_1\cdots x_{d+1}$. The map $\phi$ is a diffeomorphism between each hyperoctant of $\R^{d+1}$ and $S$, meaning that $X_1$ restricted to each hyperoctant is equivalent to $X_0$ restricted to $S$. Therefore, $X_1$ has at least $\tau_1=2^{d+1}\cdot \tau_0$ normally hyperbolic invariant tori, $\tau_0$ of them in each of the $2^{d+1}$ hyperoctants of $\R^{d+1}$. This process can be continued to create a sequence of polynomial vector fields $(X_k)$ of degree $m_k$ with at least $\tau_k$ normally hyperbolic invariant tori, where $m_k=2 m_{k-1}+d$ and $\tau_k=2^{d+1} \tau_{k-1}$. This yields:
\[
m_k=2^{k}(m_0+d)-d\quad\text{and}\quad \tau_k=2^{k(1+d)}\tau_0=\dfrac{\tau_0}{(d+m_0)^{d+1}}(m_k+d)^{d+1}.
\]
Thus, 
\[
N_h^d(m_k)\geq \dfrac{\tau_0}{(d+m_0)^{d+1}}(m_k+d)^{d+1}>\dfrac{\tau_0}{(d+m_0)^{d+1}}m_k^{d+1},
\]
implying that 
\begin{equation}\label{liminf2}
\lim_{m\to\infty}\inf\dfrac{N_h^d(m)}{m^{d+1}}\geq\dfrac{\tau_0}{(d+m_0)^{d+1}}.
\end{equation}
Notice that, if the $\tau_0$ invariant $d$-tori of the initial polynomial vector field $X_0$ are considered isolated but not necessarily normally hyperbolic, the relationship \eqref{liminf2} would still hold, but for $N^d(m)$ instead.

Therefore, once we obtain an example of a $(d+1)$-dimensional polynomial vector field $X_0$ of degree $m_0$ with at least $\tau_0$ normally hyperbolic (resp. isolated) invariant $d$-tori, the relationship \eqref{liminf2} establishes that $N_h^d(m)$ (resp. $N^d(m)$) grows at least as fast as $(\tau_0/(d+m_0)^{d+1})\,m^{d+1}$. To obtain such an example, one possible approach would be to extend Theorem \ref{thm:nhit} to detect higher-dimensional tori.

\bibliographystyle{abbrv}
\bibliography{references.bib}
\end{document}